\documentclass[11pt]{amsart}
\usepackage{epsfig}
\usepackage{verbatim}
\usepackage{amssymb}
\usepackage{amsfonts}
\usepackage{txfonts}
\usepackage{graphicx}
\newcommand{\avg}{\mathrm{avg}}
\newcommand{\var}{\mathrm{var}}
\newcommand{\clus}{\mathrm{clus}}

\newcommand{\bvec}{\left[ \begin{array}{c}}
\newcommand{\evec}{\end{array} \right]}
\newcommand{\refl}{\mathrm{refl}}

\newcommand{\ra}{{\rightarrow}}

\newcommand{\eproof}{\hfill\rule{2.2mm}{3.0mm}}
\newcommand{\esubproof}{\hfill$\Box$}
\newcommand{\Proof}{\noindent {\bf Proof.~~}}

\newcommand{\HH}{{\mathcal H}}

\renewcommand{\SS}{{\mathcal S}}

\newcommand{\Z}{{\mathbb Z}}

\renewcommand{\eqref}[1]{(\ref{#1})}

\newcommand{\wgt}{\mathrm{wgt}}

\newcommand{\bfs}{{\mathbf s}}

\newcommand{\bft}{\mathbf t}
\newtheorem{prop}{Proposition}[section]
\newtheorem{lem}[prop]{Lemma}

\newtheorem{coro}[prop]{Corollary}
\newtheorem{theo}[prop]{Theorem}
\newtheorem{conj}[prop]{Conjecture}
\newtheorem{exam}{Example}[section]

\begin{document}
\baselineskip 18pt
\title{The Structure of a Bernoulli Process Variation of the Fibonacci Sequence}
\author{Brian A. Benson}

\email{gth858n@mail.gatech.edu}

$
$
\begin{abstract}
We consider the structure of a variation of the Fibonacci sequence which is determined by a Bernoulli process.  The associated structure of all Bernoulli variations of the Fibonacci sequence can be represented by a directed binary tree, which we denote $X$, with vertex labels representing the specific state of the recurrence variation.  Since $X$ is a binary tree, we can consider the term of a sequence variation given by a finite traversal of $X$ represented by a binary code $\bft$.  We then prove that the traversal of $X$ that is the reflection of the digits of $\bft$ gives exactly the integer term corresponding to $\bft$.  We consider how to further this result with the statement of an additional conjecture.  Finally, we give connections to Fibonacci expansions, the Stern-Brocot tree, and we apply our methods to the Three Hat Problem as seen in {\em Puzzle Corner} of the {\em Technology Review} magazine.

\end{abstract}

\maketitle

\section{Introduction}
\setcounter{equation}{0}

Variations of the Fibonacci sequence arise in several applied fields of study\cite{BK06}, for example phyllotaxis \cite{J84, J94}.  In \cite{BK06}, subtle variations of the Fibonacci sequence are explored from a probabalistic perspective.  Herein, we consider a particular variation of the Fibonacci sequence determined by a Bernoulli process.  First, we describe the combinatorial structure of these variations using binary code associated with the independent random variables of the Bernoulli process.  Specifically, we prove that for any such code of finite length, then the digits of the code listed in reverse order correspond to the same integer term in a distinct variation of the Fibonacci sequence.  We then give a conjecture as to the partial ordering of the terms corresponding to codes of the same length and weight.

To further motive the work herein, we give several connections corresponding to the structure of the Bernoulli process variation of the Fibonacci sequence.  In fact, the structure created by the Bernoulli process variation on the Fibonacci sequence can be associated with a specific collection of Fibonacci expansions in a very natural way; specifically, we give a simple correspondence between the structure of Bernoulli variations of the Fibonacci sequence and this collection of Fibonacci expansions.   In addition, we give a connection between the structure of this variation and the Stern-Brocot tree.  Finally, the variations of the Fibonacci sequence given herein can be applied to the Three Hat Problem as seen in {\em Puzzle Corner} of the {\em Technology Review} magazine.

\section{A Bernoulli Variation of the Fibonacci Sequence}
\setcounter{equation}{0}

In this section, we construct a sequence $(B_m)$ which is a variation of the Fibonacci sequence $(F_m)$, by an associated Bernoulli process; we refer to the following construction as a {\em Bernoulli variation} of the Fibonacci sequence.  We first consider the construction of the Bernoulli process.  Let $X=\{0, 1 \}^{\Z^+}$ and if $\SS$ is the $\sigma$-algebra of the finite set $\{0, 1 \}$, then let $\Sigma = \SS^{\Z^+}$. Let $p_i$ be the probability associated with some $i \in \{0, 1 \}$ such that $p_0 + p_1 = 1$.  Letting $\mu = \{p_0 , p_1 \}^{\Z^+}$ then $(X , \Sigma , \mu)$ comprises the probability space defining our Bernoulli process\footnote{Note that we will not consider any results with respect to probabilities herein.}.

To construct the sequence $(B_m)$, first, let $B_k = F_k$ for integers $k = 1,2,3$.  Let the {\em root state} $S_1$ be the point $[F_1, F_2, F_3]$ of $\Z^{n+1}$.  The variation of the recurrence scheme generated by the Bernoulli scheme is given by the following construction:  for $(x_1, x_2, \ldots , x_k , \ldots) \in X$ and for integer $k \geq 2$, the state $S_k$ is given recursively by $$S_k \coloneqq \tau_{x_{k-1}}(S_{k-1})$$  where $\tau_i : \Z^{n+1} \ra \Z^{n+1}$ is is defined by \vspace{3mm}\begin{small}$$ \tau_i ([B_{k+1}, B_{k+2}, B_{k+3}])= \left \{ \begin{array}{lr} 
	\left[ B_{k+1}, B_{k+3}, B_{k+1}+B_{k+3} \right ] , & \mbox{if $i=0$;} \\ 
	 \left [ B_{k+2}, B_{k+3}, B_{k+2} + B_{k+3} \right ], & \mbox{if $i=1$;}\\ \end{array} \right. $$ \end{small}  Although we do not take full advantage of the probablistic potential of our construction herein, due to the number of connections associated with this variation, authors of future work related to these connections may wish to consider such approaches.
	 
As a final note on this construction, when we wish to consider all Bernoulli variations of the Fibonacci sequence, we consider a directed, binary tree where each vertex represents a state of all possible Bernoulli variations of $(F_m)$.  This tree is a representation of $X$ where the elements of $X$ are infinite traversals of the directed edges of the tree; more specifically, $(x_1, x_2, \ldots, x_k , \ldots) \in X$ represents an infinite traversal of the edges $x_1, x_2, \ldots, x_k, \ldots$ of the tree.  The root of the tree is $S_1$ and the out-degree of each vertex is $n$ corresponding to the $n$ possible values of each $x_k$; thus, we can label each of the out edges of a vertex with the appropriate $x_k$.  Henceforth, we will often refer to $X$ as this tree.

\begin{figure}
\begin{center}
	\includegraphics[width=3.5in]{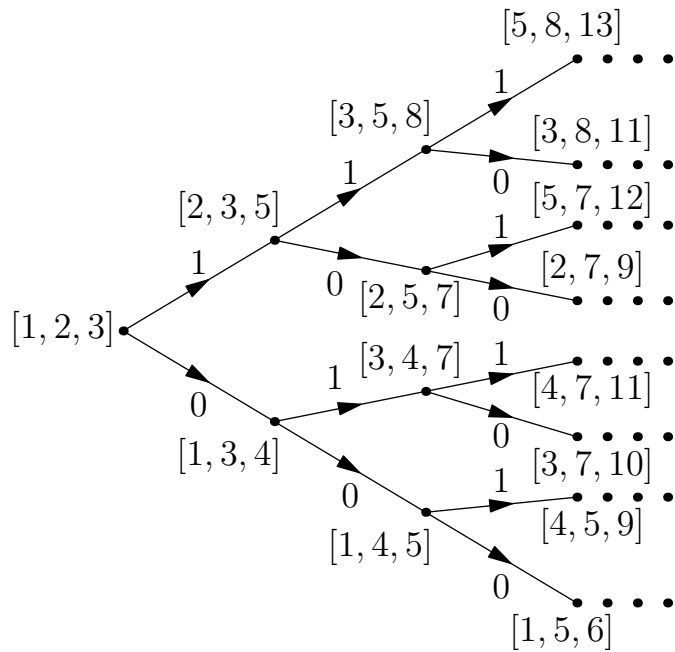}
\end{center}
\label{F1}
\caption{The State Representation of the Bernoulli Variation of the Fibonacci Sequence, $X$}
\end{figure}

\begin{prop} The map $\tau_i$ on is additively linear on vectors of the form $[a,b,a+b]$ with $a,b \in \Z^+$.\end{prop}

\Proof Let $a,b,a',b' \in \Z^+$.  For $i = 0$, we have $$\tau_0([a, b, a+b])+\tau_0([a', b', a'+b'])=[a, a+b, 2a+b] + [a', a'+b', 2a' +b']=$$ $$[a+a', a+b+a'+b', 2a+b+2a'+b']= \tau_0 ([a+a', b+b', a+a' + b +b'])=$$ $$\tau_0 ([a,b,a+b]+[a', b', a'+b']).$$ If $i=1$, we have $$\tau_1([a,b,a+b])+\tau_1([a', b', a'+b']) = [b,a+b,a+2b]+[b', a'+b', a'+2b'] =$$ $$ [b+b', a+b+a'+b', a+2b+a'+2b']=\tau_1([a+a', b+b', a+a'+b+b'])=\tau_1([a, b, a+b]+[a',b',a'+b']).$$\eproof

\noindent {\bf Notation.}  We refer to $B_{n+k+1}$ corresponding to $(x_1, x_2, \ldots , x_k , \ldots ) \in X$ as $F[x_1x_2 \cdots x_k]$\footnote{Note that by construction, if $x_i = n-1$ for all $i \in \Z^+$, then $(B_m)=(a_m)$.}; this is so we can reference all possible Bernoulli variations of $(a_m)$ in .  Similarly, we can consider the entire state corresponding to $F[x_1x_2 \cdots x_k]$ which we denote $S[x_1x_2 \cdots x_k]$.

\begin{prop}
Let $A=\{\bfs \in \Z^{3}: \bfs = [a,b,a+b], a, b \in \Z_{\geq 0}\}$.   For each $S_k$, $k \geq 2$, there exists a function $r:A \ra \Z^{3}$ such that $r(S_k)=S_{k-1}$.  Further, $$r([a,b,c])= \sigma\left ( [a,b-a,b] \right )$$ for some $\sigma$, a dimensional permutation of $\Z^{3}$ which fixes the third entry of the vector.
\end{prop}

\Proof By construction, we know that $a$ and $b$ are entries of the vector $S_{k-1}$.  Further, by construction, we know that $b$ is the final entry of $S_{k-1}$.  Thus, there is one entry of $S_{k-1}$ which is unaccounted for.  Let $y$ be this entry, then $y+a =b$ which tells that $y = b-a$.  Thus, $S_{k-1} = \sigma [a, b-a, b ]$ for some dimension permutation $\sigma$ of $\Z^{3}$ which fixes the third entry of $[a, b-a, b]$.\eproof

We again denote the sample space of the Bernoulli variation of the Fibonacci sequence as $X$ and denote the element of the Bernoulli variation of the Fibonacci sequence given by the traversal $x_1 x_2 \cdots x_k$ of $X$ as $F[x_1 x_2 \cdots x_k]$.  Further, for shorthand, we denote a traversal of $X$ as $\bft =x_1x_2 \cdots x_k$.  We define the {\em reflection} of a traversal $\bft=x_1x_2 \cdots x_k$ denoted $\refl (\bft)$ or $\refl(x_1x_2 \cdots x_k)$ as the permutation on that traversal corresponding to $x_kx_{k-1}\cdots x_1$; more directly stated, 
$$\refl(x_1x_2 \cdots x_k)=x_kx_{k-1}\cdots x_1.$$  This idea leads to one of our main results which tells us that $F[\bft]=F[\refl (\bft)]$.  However, before we explicitly state and prove this result, we illustrate the basic idea with a few examples.

\begin{exam}$F[1011]=F[1101]$\end{exam}
To show this example, we give the computations of $S[1011]$ and $S[1101]$ respectively.  First, $S[1011]$ has the traversal $$[1,2,3]\stackrel{\tau_1}{\ra}[2,3,5]\stackrel{\tau_0}{\ra}[2,5,7]\stackrel{\tau_1}{\ra}[5,7,12]\stackrel{\tau_1}{\ra}[7,12,19]$$ implying that $F[1011]=19$.  Second, $S[1101]$ has the traversal $$[1,2,3]\stackrel{\tau_1}{\ra}[2,3,5]\stackrel{\tau_1}{\ra}[3,5,8]\stackrel{\tau_0}{\ra}[3,8,11]\stackrel{\tau_1}{\ra}[8,11,19]$$ implying that $F[1101]=19$.\esubproof

\begin{exam} $F[1010000]=F[0000101]$\end{exam}
Note that $S[1010000]$ corresponds to the traversal
\begin{small}$$[1,2,3]\stackrel{\tau_1}{\ra}[2,3,5]\stackrel{\tau_0}{\ra}[2,5,7]\stackrel{\tau_1}{\ra}[5,7,12] \stackrel{\tau_0}{\ra}[5,12,17]\stackrel{\tau_0}{\ra}[5,17,22]\stackrel{\tau_0}{\ra}[5,22,27]\stackrel{\tau_0}{\ra}[5,27,32]$$\end{small}
while $S[0000101]$ corresponds to the traversal
\begin{small}$$[1,2,3]\stackrel{\tau_0}{\ra}[1,3,4]\stackrel{\tau_0}{\ra}[1,4,5]\stackrel{\tau_0}{\ra}[1,5,6]\stackrel{\tau_0}{\ra}[1,6,7]\stackrel{\tau_1}{\ra}[6,7,13]\stackrel{\tau_0}{\ra}[6,13,19]\stackrel{\tau_1}{\ra}[13,19,32].$$\end{small}  Thus, $F[1010000]=32=F[0000101]$.\esubproof

\begin{theo}(Traversal Reflection)
If $\bft$ is a code associated with a finite traversal of $X$, then $F[\bft] = F[\refl(\bft)]$.
\end{theo}

\Proof  Clearly, the theorem is true if a traversal is a palindrome; thereby, we can assume that all traversals henceforth are not palindromes.

For a Fibonacci state $[a,b,c]$ with $a<b<c$, note that by proposition 2.2, we have $r([a,b,c])=\sigma[a,b-a,b]$ where we can let $\sigma$ be a well-ordering permutation of dimensions as a result of the fact that all Bernoulli variations of $(F_m)$ are monotone.  In other words, for the Fibonacci sequence, $\sigma$ maintains the ordering of the state under the well-ordering of the entries.  More specifically, $\sigma$ ensures that $r([a,b,c])=[a,b-a,b]$ when $a \leq b-a$ and $r([a,b,c])=[b-a,a,b]$ when $b-a \leq a$.

Note that by proposition 2.1, $\tau_0$ and $\tau_1$ are linear since the function which generates the Fibonacci sequence is linear.  Note that figure 1 above provides us with a base case for induction since it gives us the cases for traversal length equal to 2.  Suppose that the reflection principle is true for all integers $k \leq n-1$. Further, consider the traversal $x_1x_2 \cdots x_n$.  To proceed, we must consider $\tau_{x_n}(S[x_1x_2 \cdots x_{n-1}])$; from here, we consider two cases, $x_n=0$ and $x_n=1$.

Due to several complications in the case where $x_n=0$, we initially suppose $x_n=1$\footnote{We will be able to use our proof of the simpler case of $x_n=1$ in order to simplify the number of cases that we must consider when $x_n=0$.}; we note that we can represent $\tau_1$ on a state $\bfs=[a,b,c]$ as $\tau_1(\bfs)=\tau_1 ([a,b,c])=[b,c,b+c]=[a,b,c]+[b-a,a,b]=\bfs+\sigma^{\ast}(r(\bfs))$ up to some permutation $\sigma^{\ast}$ which fixes the third entry of the vector\footnote{Note that while $\sigma{\bfs}$ arranges the entries of $\bfs$ by well-ordering the integers from least to greatest, $\sigma^{\ast}(\bfs)$ only requires that we fix the third entry.  Thus, in general, $\sigma (\bfs ) \neq \sigma^{\ast} (\bfs)$.}.  Therefore, $$S[x_1x_2 \cdots x_n]=\tau_1(S[x_1x_2\cdots x_{n-1}])=S[x_1x_2\cdots x_{n-1}]+\sigma^{\ast}(r(S[x_1x_2\cdots x_{n-1}]))=$$ $$ S[x_1x_2\cdots x_{n-1}]+\sigma^{\ast}(S[x_1x_2\cdots x_{n-2}]).$$  Now, $$S[x_nx_{n-1}\cdots x_2x_1]=\tau_{x_1}\circ \tau_{x_2} \circ \cdots \circ \tau_{x_{n-1}} \circ \tau_{x_n} \left(\bvec 1\\2\\3\\ \evec \right)=$$ $$\tau_{x_1}\circ \tau_{x_2} \circ \cdots \circ \tau_{x_{n-1}}\left(\bvec 2\\3\\5\\ \evec \right)=\tau_{x_1}\circ \tau_{x_2} \circ \cdots \circ \tau_{x_{n-1}}\left( \bvec 1\\1\\2\\ \evec + \bvec 1\\2\\3\\ \evec \right) =$$ $$\tau_{x_1}\circ \tau_{x_2} \circ \cdots \circ \tau_{x_{n-1}}\left(\bvec 1\\1\\2\\ \evec \right) + \tau_{x_1}\circ \tau_{x_2} \circ \cdots \circ \tau_{x_{n-1}}\left(\bvec 1\\2\\3\\ \evec \right)=$$ $$S[x_{n-2}x_{n-3} \cdots x_2 x_1]+S[x_{n-1}x_{n-2} \cdots x_2 x_1].$$  By the induction hypothesis, $F[x_{n-2}\cdots x_2x_1]=F[x_1x_2 \cdots x_{n-2}]$ and $F[x_{n-1}\cdots x_2x_1]=$

\noindent
$F[x_1x_2 \cdots x_{n-1}]$.  Since each of these values are in the third vector position of the equations above\footnote{Note that our permutation $\sigma^{\ast}$ fixed the third entry of a vector.}, we have that $$F[x_1x_2 \cdots x_n]=F[x_1x_2 \cdots x_{n-2}]+F[x_1x_2 \cdots x_{n-1}]=$$ $$F[x_{n-2}\cdots x_2x_1]+F[x_{n-1}\cdots x_2x_1]=F[x_{n}\cdots x_2x_1]$$ proving the case of $x_n=1$.

Now, suppose $x_n=0$\footnote{One of the main difficulties with this case arises from the fact that $\tau_0[0,1,1]=[0,1,1]$, so we must alter our argument from the case where $x_n=1$.}, clearly, if $x_1=1$, then the proof of the case where $x_n=1$ will suffice in proving this case as well.  Thus, we can assume that $x_1,x_n=0$.  Without loss of generality, we can assume that $x_{n-k}=1$, $1 \leq k \leq \lfloor n/2 \rfloor$, and for all $j < k$, $x_{n-j} =0$ and $x_j =0$.  To prove the case of $x_n=0$, we will induct on $k$ beginning with the base case of $k=1$.

When $k=1$, we know that $x_1=x_n=0$ and $x_{n-1}=1$.  Then, if $S[x_1x_2 \cdots x_{n-1}]=[a,b,c]$, then $$S[x_1x_2 \cdots x_n] = S[x_1x_2 \cdots x_{n-1}] + \bvec 0\\a\\a\\ \evec.$$  Further, $$S[x_nx_{n-1} \cdots x_1]=\tau_{x_1} \circ \tau_{x_2} \circ \cdots \circ \tau_{x_{n}} \left( \bvec 1\\2\\3\\ \evec \right)=$$ $$\tau_{x_1} \circ \tau_{x_2} \circ \cdots \circ \tau_{x_{n-1}}\left( \bvec 1\\2\\3\\ \evec + \bvec 0\\1\\1\\ \evec \right)=S[x_{n-1}x_{n-2} \cdots x_1]+\tau_{x_1} \circ \cdots \circ \tau_{x_{n-1}}\left( \bvec 0\\1\\1\\ \evec \right).$$

\vspace{2mm}
\noindent
Since $\tau_{x_{n-1}}=1$, we have $\tau_{x_1} \circ \cdots \circ \tau_{x_{n-1}}([0,1,1])=\tau_{x_1} \circ \cdots \circ \tau_{x_{n-2}}([1,1,2])=\tau_{x_1} \circ \cdots \circ \tau_{x_{n-3}}([1,2,3])=S[x_{n-3}x_{n-4} \cdots x_1]$.  Thus, we have $$S[x_nx_{n-1} \cdots x_1]=S[x_{n-1}x_{n-2} \cdots x_1]+S[x_{n-3}x_{n-4} \cdots x_1].$$

Since, by the original induction hypothesis, $F[x_1x_2 \cdots x_{n-1}]=F[x_{n-1}x_{n-2} \cdots x_1]$, we must check that $F[x_{n-3}x_{n-4} \cdots x_1]=a$.  To do this, we again rely on the original induction hypothesis to tell us that $F[x_{n-3}x_{n-4} \cdots x_1]=F[x_1x_2 \cdots x_{n-3}]$\footnote{This allows us to compute the final entry of $S[x_1x_2 \cdots x_{n-3}]$ in place of the final entry of $S[x_{n-3}x_{n-4} \cdots x_1]$.}.

Now, since we took $S[x_1x_2 \cdots x_{n-1}]=[a,b,c]$, we have that $$S[x_1 \cdots x_{n-2}]=r(S[x_1 \cdots x_{n-1}]).$$  Now, $r(S[x_1 \cdots x_{n-1}])= [a,b-a,b] \text{ or }[b-a,a,b]$ depending on the ordering of $a$ and $b-a$.  However, since we assumed that $x_{n-1}=1$, we consider the equation $$\tau_{x_{n-1}} \left( \bvec a'\\b'\\c'\\ \evec \right) = \tau_{1} \left( \bvec a'\\b'\\c'\\ \evec \right)=\bvec b'\\c'\\b'+c'\\ \evec = \bvec a\\b\\c\\ \evec$$ where $a'+b'=c'$ and $a'<b'<c'$.  By our equation, we have that $a=b', b=c', c=b'+c'$ which implies the equalities $a'=c'-b'=b-a, b'=a, c'=c-a=b$.  Therefore, $b-a <a<b$ implying that $$S[x_1 \cdots x_{n-2}]=r(S[x_1 \cdots x_{n-1}])= r\left( \bvec a\\b\\c\\ \evec \right)= \bvec b-a\\a\\b\\ \evec.$$  Now, $S[x_1 \cdots x_{n-3}]=r(S[x_1 \cdots x_{n-2}])=r[b-a,a,b]=[|b-2a|,b-a,a]$ or $[b-a, |b-2a|,a]$ depending on the ordering of $b-a$ and $|b-2a|$.  Either way, we can conclude that $F[x_{n-3} \cdots x_1] = F[x_1 \cdots x_{n-3}] = a$.  Thus, we can conclude when $x_n=0$ and $k=1$, $F[x_1 \cdots x_n]=F[x_n \cdots x_1]$ completing the base case of the induction on $k$.

Now, assume that the claim is true for all integers $k$ such that $1<k<m< \lfloor n/2 \rfloor$.  Now, we consider the case of $m+1$.  Then, by this, we know that $x_1, \ldots, x_{m-1}=0$ and $x_{n-m}, \ldots, x_{n}=0$.  If we again take $S[x_1 \cdots x_{n-1}]=[a,b,c]$, we have that $$S[x_1 \cdots x_n]=S[x_1 \cdots x_{n-1}] + \bvec 0\\a\\a\\ \evec.$$  Further, $$S[x_n\cdots x_1]=S[x_{n-1} \cdots x_1]+\tau_{x_1} \circ \cdots \circ \tau_{x_{n-1}}\left( \bvec 0\\1\\1\\ \evec \right).$$  However, since $x_{n-m}, \ldots, x_{n-1}=0$, $\tau_{x_1} \circ \cdots \circ \tau_{x_{n-1}}([0,1,1])=\tau_{x_1} \circ \cdots \circ \tau_{x_{n-m+1}}([0,1,1])=S[x_{n-m-1} \cdots x_1]$.  Similar to the base case, we consider the third entry of $S[x_1 \cdots x_{n-m-1}]=$

\noindent
$r^m(S[x_1 \cdots x_{n-1}])$ which is equal to $F[x_{n-m-1}\cdots x_1]$ under the induction hypothesis.

Again refering to the configuration $[a',b',c']$ with $a'<b'<c'$, we have the we wish to find the exact permutation $\sigma$ required for $r([a,b,c])$ where $\tau_0([a',b',c'])=[a,b,c]$.  Thus, we consult the equation $$\tau_0 \left( \bvec a'\\b'\\c'\\ \evec \right) = \bvec a'\\b'\\c'\\ \evec + \bvec 0\\a'\\a'\\ \evec = \bvec a\\b\\c\\ \evec.$$  Solving for $a'$, $b'$, and $c'$, we get $a'=a$, $b'=b-a$, and $c'=c-a$.  Since $x_{n-m}, \ldots, x_{n}=0$, we know that $\min(S[x_1 \cdots x_{n-m+1}]=\min(r^{m-2}(S[x_1 \cdots x_{n-1}]))=a$.  To find the third entry of $S[x_1 \cdots x_{n-m-1}]$, we must find $\max(r^2(S[x_1 \cdots x_{n-m+1}))$; but since $x_{n-m+1}=1$, we have that $\max(r^2(S[x_1 \cdots x_{n-m+1}))=\min(S[x_1 \cdots x_{n-m+1}]=a$ by our proof of the base case of the induction on $k$.  Thus, $F[x_1 \cdots x_n]=F[x_1 \cdots x_{n-1}]+a = F[x_{n-1} \cdots x_{1}]+a=F[x_n \cdots x_1]$ completing the proof of the case $x_n=0$.  Thus, all possible cases of non-palindrome traversals have been exhausted proving the theorem.\eproof

\noindent
{\bf Remark.}  The converse of the reflection principle, for traversals of equal length, does not hold true in general.  Consider the traversals $10011$ and $01110$; a simple  compution yeilds $F[10011]=25=F[01110]$ confirming this fact.

In general, suppose that we consider Bernoulli variations generated by an arbitrary state $[a,b,c]$ where $a,b,c \in \Z^+$ and $a+b=c$.  The following result tells us that the primitive configuration $[1,2,3]=[F_2, F_3, F_4]$ holds special significance with respect to the reflection principle.

\begin{prop}
The only two pair-wise primitive root configurations in which the reflection theorem holds are $[1,2,3]$ and $[2,1,3]$.
\end{prop}

\Proof If the reflection theorem is true for an arbitrary, pair-wise primitive root configuration $[a,b,c]$ such that $a,b,c \in \Z^+$ with $a+b=c$, then we must have that $F[01]=F[10]$.  First assume that $a\leq b\leq c$.  Since $S[01]=[c,a+c,a+2c]$ and $S[10]=[b,b+c,2b+c]$, we have that $F[01]=a+2c$ and $F[10]=2b+c$.  Therefore, for the reflection principle to hold, we must have $a+2c=2b+c$.  Since $a+b=c$, we have that $2a=b$.  Since $[a,b,c]$ is pair-wise primitive, we must have that $a=1$ and $b=2$ giving us the result for $[1,2,3]$.  Similarly, if we assume that $b\leq a \leq c$, we get that $a=2b$, so by pair-wise primitivity, we have that $a=2$ and $b=1$ giving us the result for $[2,1,3]$.  Note that because $G_1, G_2 \in \Z^+$, this exhausts all other possibilities.  Finally, note that in the case of $[2,1,3]$ as the root, then any code $x_1 x_2\cdots$ in $X$ corresponds to the code $y_1y_2 \cdots $ in $X$ where $y_i = 1-x_i$; thus, the reflection property still holds giving the result.\eproof

Motivated by the reflection theorem, we would like to see how the number of changes or transitions of a code $x_1x_2 \cdots x_k$ affects the magnitude of the integer associated with $F[x_1 x_2 \cdots x_k]$.

\begin{prop}
For integer $j \geq 2$, $$F[\underbrace{00\cdots 0}_{j}\underbrace{11\cdots 1}_{j}]=F[\underbrace{11\cdots 1}_{j}\underbrace{00\cdots 0}_{j}]<F[\underbrace{0101 \cdots 01}_{2j}]=F[\underbrace{1010 \cdots 10}_{2j}].$$
\end{prop}

\Proof The equivalent conditions hold by theorem 3.1; further, it is only necessary to compare $F[\underbrace{0101\cdots01}_{2j}]$ and $F[\underbrace{11\cdots 1}_{j}\underbrace{00\cdots 0}_{j}]$. The base case is $j=2$ where $S[1100]=[3,11,14]$ and $S[0101]=[7,10,17]$; from here we note that $\min(S[1100])<\min(S[0101])$ and $\max(S[1100])<\max(S[0101])$.  Now, we apply the step of induction assuming that these minimum and maximum inequality statements are true for all integers $j\leq k$ for arbitrary $k\in \Z^+$.  Since $S[1]=[F_3,F_4,F_5]$, where $F_i$ corresponds to the $i$-th integer in the Fibonacci sequence, we note that $S[\underbrace{11\cdots 1}_{j}\underbrace{00\cdots 0}_{j}]=[F_{j+2},F_{j+4}+(j-1)F_{j+2},F_{j+4}+jF_{j+2}]$.

Let $S[\underbrace{0101\cdots01}_{2k}]=[a,b,c]$\footnote{Note that we are referring to a general configuration.} with $a+b=c$.  Note that $S[\underbrace{0101\cdots01}_{2(k+1)}]=[c,a+c,a+2c]$ and $S[\underbrace{11\cdots 1}_{k+1}\underbrace{00\cdots 0}_{k+1}]=[F_{k+3},F_{k+3}+kF_{k+3},F_{k+5}+(k+1)F_{k+3}]$.  By the induction hypothesis and the fact that the Fibonacci sequence is monotone increasing, we have $a>F_{k+2}>F_{k+1}$.  Thus, $$a+2c=a+2(a+b)>a+2(2a)=5a>2a>2F_{k+2}>F_{k+2}+F_{k+1}=F_{k+3}.$$  By this, we have $$F_{k+5}+(k+1)F_{k+3}=F_{k+4}+F_{k+3}+(k+1)(F_{k+2}+F_{k+1})=$$ $$F_{k+4}+kF_{k+2}+F_{k+3}+F_{k+2}+(k+1)F_{k+1}<c+F_{k+3}+F_{k+2}+(k+1)F_{k+1}=$$ $$c+F_{k+4}+kF_{k+1}+F_{k+1}<c+F_{k+4}+kF_{k+2}+F_{k+2}<2c+a.$$ This gives the result.\eproof

To give a more generalized notion of this observation, we define several properties of general traversals of $X$. Consistent with binary codes in coding theory, we let the weight of a configuration, $\wgt (\bft)$, represent the number of ones present in the traversal.  We define the {\em edge cluster number} of the edge $x_i$ in the traversal $\bft=x_1x_2\cdots x_i \cdots x_n$ to be $\clus (x_i, \bft)=|\{x_k:x_k=x_j=x_i \text{ for all } \min(i,k) \leq j \leq \max(i,k)\}|$.  From this, we can define the {\em cluster average} of the traversal $\bft$ of $X$ to be $$\avg(\bft)=\sum_{i=1}^{n}{\frac{[\clus(x_i,\bft)]}{n}}.$$  Further, we define the {\em cluster variance} of the traversal $\bft$ of $X$ to be $$\var(\bft)=\sum_{i=1}^{n}{\frac{[\clus(x_i,\bft)]^2}{n}}.$$

The following is an example of the relation between the variance of a traversal code and its correponding value with respect to the Bernoulli variation of Fibonacci sequence.

\begin{exam} Note that $\wgt(\underbrace{0101\cdots 01}_{2j})=\wgt(\underbrace{11\cdots 1}_{j}\underbrace{00\cdots 0}_{j})$ while $\var(\underbrace{0101\cdots01}_{2j})=2j<2j^3=j\cdot j^2+j\cdot j^2=\var(\underbrace{11\cdots 1}_{j}\underbrace{00\cdots 0}_{j})$.  Further, by the proposition, $F[\underbrace{11\cdots 1}_{j}\underbrace{00\cdots 0}_{j}]< F[\underbrace{0101\cdots 01}_{2j}]$.\esubproof \end{exam}
\noindent
This example illustrates the beginning of observations which give evidence for the following conjecture.  The remaining rationale for the conjecture is given after its statement.

\begin{conj}
For traversals $\bft_1,\bft_2$ of $X$ of equal length, if $\wgt(\bft_1)=\wgt(\bft_2)$ and $\var(\bft_1)<\var(\bft_2)$, then $F[\bft_1]>F[\bft_2]$.
\end{conj}

The overriding rationale behind the conjecture is that the larger the variance of the code of a traversal, the smaller the average clusters size and, thus, the more transitions there are back and forth between ones and zeros when the weight or the code is constant.  Now, suppose that we wish to maximize $F[\bft]$ for traversal $\bft$ of fixed length and weight; then, we suppose that $S[x_1 x_2 \cdots x_k]=[a,b,c]$ which naturally implies that $a,b \in \Z^+$, $a < b$, and $a+b=c$.  Now note that $S[x_1 x_2 \cdots x_k 0]=[a,c,a+c]$ while $S[x_1 x_2 \cdots x_k1]=[b,c,b+c]$ which means that the first entry of $S[x_1 x_2 \cdots x_kx_{k+1}]$ is maximized locally by choosing $x_{k+1}=1$ while the second entry of $S[x_1 x_2 \cdots x_kx_{k+1}]$ is the same irregardless of the value of $x_{k+1}$.  Since the length and weight of $\bft$ are fixed, when $x_j$ must be zero, having $x_{j-1}=1$ maximizes the sequence locally.  Spreading this local observation over the entire length of the code gives evidence for the conjecture.  Further, although these local observations are relatively straightforward, it appears that constructing a rigorous proof of the conjecture from these observations is somewhat less intuitive.  In the following paragraphs of this section, we will consider a few approaches towards proving the conjecture.

Perhaps the first approach a reader might take is straightforward induction on the code length\footnote{In addition, perhaps even induction on the code weight for each code length as well.}.  However, if $\wgt(x_1 \cdots x_k)=\wgt(x_1' \cdots x_k')$ and $\var [x_1 \cdots x_k] < \var [x_1' \cdots x_k']$, then it is not necessarily true that $\var [x_1 \cdots x_k x_{k+1}] \leq \var [x_1' \cdots x_k'x_{k+1}]$\footnote{Note, however, that if it is the case that $x_k \neq x_{k+1}$, then $\var [x_1 \cdots x_k x_{k+1}]=(n/ (n+1))(\var [x_1 \cdots x_k]+1) < (n/(n+1))(\var [x_1' \cdots x_k']+1) \leq \var [x_1' \cdots x_k'x_{k+1}]$.}.  A counter-example which tells us that this is not true in general is as follows:  $\var[1010111] = 31/7 < 55/7 =\var [1110110]$, however, $\var [10101111] = 17/2 > 65/8 = \var[11101101]$.  

Note that even if the orderings of the variance of $\bft$ is known with respect to an inductive step, this in itself does not appear to be enough to establish the ratio between the first and second entry of $S[\bft]$; to determine this ratio from a label without direct computation from the code, it appears that something else must be known about the structure of the code.  However, as the length of the codes under consideration becomes larger, additional structures arise which make this approach non-trivial.

\begin{exam} We consider all traversal codes of length 4 within the context of the conjecture.\end{exam}

The traversal codes of the non-trivial weights are given below.

\begin{center}
\begin{tabular}[t]{ccccc}
 Traversal Code $\bft$ & $\var [\bft]$ & Generation of $S[\bft]$ \\
 $1000$ & 7 &  $[1,2,3] \stackrel{\tau_1}{\rightarrow} [2,3,5] \stackrel{\tau_0}{\rightarrow} [2,5,7] \stackrel{\tau_0}{\rightarrow} [2,7,9] \stackrel{\tau_0}{\rightarrow} [2,9,11]$\\
 $0100$ & 5/2 & $[1,2,3] \stackrel{\tau_0}{\rightarrow} [1,3,4] \stackrel{\tau_1}{\rightarrow} [3,4,7] \stackrel{\tau_0}{\rightarrow} [3,7,10] \stackrel{\tau_0}{\rightarrow} [3,10, 13]$\\
 $0010$ & 5/2 & $[1,2,3] \stackrel{\tau_0}{\rightarrow} [1,3,4] \stackrel{\tau_0}{\rightarrow} [1,4,5] \stackrel{\tau_1}{\rightarrow} [4,5,9] \stackrel{\tau_0}{\rightarrow} [4,9,13]$\\
 $0001$ & 7 & $[1,2,3] \stackrel{\tau_0}{\rightarrow} [1,3,4] \stackrel{\tau_0}{\rightarrow} [1,4,5] \stackrel{\tau_0}{\rightarrow} [1,5,6] \stackrel{\tau_1}{\rightarrow} [5,6,11]$\\ 
$1100$ & 4 & $[1,2,3] \stackrel{\tau_1}{\rightarrow} [2,3,5] \stackrel{\tau_1}{\rightarrow} [3,5,8] \stackrel{\tau_0}{\rightarrow} [3,8,11] \stackrel{\tau_0}{\rightarrow} [3,11,14]$\\
$0110$ & 5/2 & $[1,2,3] \stackrel{\tau_0}{\rightarrow} [1,3,4] \stackrel{\tau_1}{\rightarrow} [3,4,7] \stackrel{\tau_1}{\rightarrow} [4,7, 11] \stackrel{\tau_0}{\rightarrow} [4,11, 15]$\\
$0011$ & 4 & $[1,2,3] \stackrel{\tau_0}{\rightarrow} [1,3,4] \stackrel{\tau_0}{\rightarrow} [1,4,5] \stackrel{\tau_1}{\rightarrow} [4,5,9] \stackrel{\tau_1}{\rightarrow} [5,9,14]$\\
$1001$ & 5/2 & $[1,2,3]\stackrel{\tau_1}{\rightarrow} [2,3,5] \stackrel{\tau_0}{\rightarrow} [2,5,7] \stackrel{\tau_0}{\rightarrow} [2,7,9] \stackrel{\tau_1}{\rightarrow} [7,9,16]$\\
$1010$ & 1 & $[1,2,3]\stackrel{\tau_1}{\rightarrow} [2,3,5] \stackrel{\tau_0}{\rightarrow} [2,5,7] \stackrel{\tau_1}{\rightarrow} [5,7,12] \stackrel{\tau_0}{\rightarrow} [5,12,17]$\\
$0101$ & 1 & $[1,2,3] \stackrel{\tau_0}{\rightarrow}[1,3,4] \stackrel{\tau_1}{\rightarrow} [3,4,7] \stackrel{\tau_0}{\rightarrow} [3,7,10] \stackrel{\tau_1}{\rightarrow} [7,10,17]$\\
$1110$ & 7 & $[1,2,3] \stackrel{\tau_1}{\rightarrow} [2,3,5] \stackrel{\tau_1}{\rightarrow} [3,5,8] \stackrel{\tau_1}{\rightarrow} [5,8,13] \stackrel{\tau_0}{\rightarrow} [5, 13, 18]$\\
$1101$ & 5/2 & $[1,2,3] \stackrel{\tau_1}{\rightarrow} [2,3,5] \stackrel{\tau_1}{\rightarrow} [3,5,8] \stackrel{\tau_0}{\rightarrow}[3,8,11] \stackrel{\tau_1}{\rightarrow} [8,11,19]$\\
$1011$ & 5/2 & $[1,2,3] \stackrel{\tau_1}{\rightarrow} [2,3,5] \stackrel{\tau_0}{\rightarrow} [2,5,7] \stackrel{\tau_1}{\rightarrow} [5,7,12] \stackrel{\tau_1}{\rightarrow} [7,12, 19]$\\
$0111$ & 7 & $[1,2,3] \stackrel{\tau_0}{\rightarrow} [1,3,4] \stackrel{\tau_1}{\rightarrow} [3,4,7] \stackrel{\tau_1}{\rightarrow} [4,7,11] \stackrel{\tau_1}{\rightarrow} [7,11,18]$\\
\end{tabular}
\end{center}

\vspace{2mm}
\noindent From this, it becomes clear that the conjecture holds true for all traversal codes of length 4; further, note that the variance alone is not enough to determine the ratio between the first and second entry of $S[\bft]$.\esubproof

\section{Connections to Fibonacci Expansions}
\setcounter{equation}{0}

We show that the variation of the Fibonacci sequence herein is directly related to a specific collection of Fibonacci expansions.  A Fibonacci expansion is a representation of an integer as the sum of Fibonacci numbers.  Although it is not related to the work herein, it is worth mentioning that one of the most well-known results with respect to Fibonacci expansions is Zeckendorf's theorem which tells us that any positive integer can be given as a unique sum of nonconsecutive Fibonacci numbers, see \cite{S52}.  For additional reading on Fibonacci expansions, see \cite{B65, H72}.

\noindent {\bf Notation:}  We denote code restrictions on a code $\bft$ where if $\bft = x_1 \cdots x_j$, then for $k$ such that $1 \leq k < j$, the restriction to the first $k$ terms is denoted $\bft |_k = x_1 \cdots x_k$.

\begin{theo} There exists a bijection between binary codes of the form $x_1 \cdots x_j$, for $j \in \Z^+$, such that $\wgt(x_1 \cdots x_j) < j$ and Fibonacci expansions of the form $aF_k +bF_{k+2}$ for all $a,b \in \Z^+$ with $\gcd(a,b)=1$, $k \in \Z_{\geq 2}$. \end{theo}

\Proof Let $\bft = x_1 \cdots x_j$.  Let $\nu(\bft)$ be the number of successive ones in $\bft$ beginnning at $x_1$; then $k = \nu (\bft) +2$ in $F_k, F_{k+2}$.  If $x_{\nu(\bft) +1} = 0$, then $a>b$ and if $x_{\nu(\bft)+1}=1$, then $b>a$.  Finally, the first two entries of $S[x_{\nu (\bft) +2} \cdots x_j]$ give the coefficients $a$ and $b$ in the ordering dictated by $x_{\nu (\bft) +1}$.

Clearly, a code of all ones gives a Fibonacci number.  At the first zero in the code $\bft$ (which, by definition, occurs at $x_{\nu(\bft)+1}$, the state becomes $[F_k, F_{k+2}, F_k + F_{k+2}]$ since $x_{\nu (\bft) +1}$ removes $F_{k+1}$ from the state in the transition from $S[\bft|_{x_{\nu(\bft)}}]$ to $S[\bft|_{x_{\nu(\bft)+1}}]$.  The entry $x_{\nu(\bft)+1}=0$ gives $S[\bft|_{x_{\nu(\bft)+1}}]=[F_k, F_k + F_{k+2}, 2F_k + F_{k+2}]$ while $x_{\nu(\bft)+1}=1$ gives $S[\bft|_{x_{\nu(\bft)+1}}]=[F_{k+2}, F_k + F_{k+2}, F_k + 2F_{k+2}]$.  Since $\tau_0, \tau_1$ are additively linear, we can view $$S[\bft]= \tau_j \circ \cdots \circ \tau_{\nu(\bft)+2} (S[\bft|_{x_{\nu(\bft)+1}}])=\left \{ \begin{array}{lr} \tau_j \circ \cdots \circ \tau_{\nu(\bft)+2} (F_k [1, 1, 2]+ F_{k+2}[0, 1, 1] ), & x_{\nu(\bft)+1}=0\\ \tau_j \circ \cdots \circ \tau_{\nu(\bft)+2} (F_{k+2}[1,1,2] + F_k [0, 1, 1]), & x_{\nu(\bft)+1}=1\\ \end{array} \right.$$

Now, suppose we have $aF_k +bF_{k+2}$.  Note that $$S[\underbrace{11 \cdots 1}_{k-2}] = [F_{k}, F_{k+1}, F_{k+2}]$$ and, thus, $$S[\underbrace{11 \cdots 1}_{k-2}0]=[F_k, F_{k+2}, F_k + F_{k+2}].$$  Further, clearly,  $$S[\underbrace{11 \cdots 1}_{k-2}0x_k]= \left \{ \begin{array}{lr} \left[ 1, 1, 2\right] F_k + \left[0, 1, 1\right] F_{k+2}, & \mathrm{ if }  x_k = 0\\ \left[1, 1, 2\right] F_{k+2}+ \left[0, 1, 1\right] F_k, & \mathrm{ if } x_k = 1\\ \end{array} \right. .$$  Since the transforms $\tau_0$ and $\tau_1$ are linear, clearly, $x_k$ determines whether $a<b$ or $a>b$ and since $F_k$ and $F_{k+2}$ are multiplied by configurations (whose first two entries are coprime), the remaining code of $x_{k+1} \cdots x_n$ clearly gives a configuration of which $a$ and $b$ are the first two entries \footnote{Note that the order of these first two entries is determined by $x_k$.}.  \eproof

\noindent {\bf Remark:} By this reasoning, the Fibonacci expansion corresponding to the code $\underbrace{11 \cdots 1}_{j-1}$ is simply $F_{j+1}+F_{j+2}=F_{j+3}$.

\begin{coro} There exists a bijection between expansions of the form $F[\bft]=F_{k_1}F_{k_2}\cdots F_{k_i}+F_{k_1\pm 2}F_{k_2 \pm 2}\cdots F_{k_i \pm 2}$ where each Fibonacci number in the expansion is greater than or equal to $F_2$  and all binary codes $\bft$ such that $\wgt(\bft) \geq 1$.\end{coro}

\Proof  Due to the striaghtforward nature of the proof, we merely give a brief overview.  Essentially, it suffices to show a bijection between $F[\bft]=F_{k_1}F_{k_2}\cdots F_{k_i}+F_{k_1\pm 2}F_{k_2 \pm 2}\cdots F_{k_i \pm 2}$ with each Fibonacci number in the expansion is greater than or equal to $F_2$ and $aF_k + b F_{k+2}$ with $gcd(a,b)=1$, $k \in \Z_{\geq 2}$.  Since $\gcd(a,b)=1$, there is a state $[a,b, a+b] \in X$ which allows us to recursively and inductively apply the theorem giving us an expansion of the form $F[\bft]=F_{k_1}F_{k_2}\cdots F_{k_i}+F_{k_1\pm 2}F_{k_2 \pm 2}\cdots F_{k_i \pm 2}$ where each Fibonacci number in the expansion is $F_2$ or greater.\eproof

\section{Connections to the Stern-Brocot Tree}
\setcounter{equation}{0}

The Stern-Brocot tree $T_{SB}$ is a binary tree which represents all nonnegative fractions.  Due to the amount of work which already exists on the Stern-Brocot tree, we give a referential overview of how $X$ is related to the Stern-Brocot tree.  The reader who wishes to learn more about the Stern-Brocot tree should see \cite{B07a, GKP90}.  The Stern-Brocot tree has a variety of simple applications in continued fractions (see \cite{B07c}).  Binary encodings of the Stern-Brocot tree are given in \cite{B07b, B07bb, GKP90}.  The first binary coding given in \cite{B07b, GKP90} directly corresponds with a canonical binary encoding of the Stern-Brocot tree in the same sense that we encoded $X$ above; that is, the encoding is with respect to a traversal beginning with the root vertex labelled $1/1$.  

To show the relation between the Stern-Brocot tree and the Bernoulli process variation of the Fibonacci sequence herein, we show a simple correspondence between the binary encoding of the Stern-Brocot tree given in \cite{B07bb} and our binary encoding of $X$.  For any state $S[\bft]=[a,b,c]$ in $X$, consider the notations $u[\bft]=a/b$ and $v[\bft]=b/a$.  Now, for any fraction on the Stern-Brocot, consider $f_R(a/b) = a/b+1 = (a+b)/b$ and $f_L(a/b) = a/(a+b)$ as defined in \cite{B07b}.  Thus, we have the correspondence such that if $u[\bft] = a/b$, then $u[\bft 0] = a/(a+b) = f_L(a/b)$ while $u[\bft 1] = b/(a+b) = 1/f_R(m/n)$.  Similarly, $v[\bft 0] = (a+b)/a = 1/f_L(a/b)$ and $v[\bft 1] = (a+b)/b = f_R(m/n)$.  Now, let $T_{SB}^{1/2}$ be the binary subtree of $T_{SB}$ rooted at the vertex $1/2$ and, similarly, let $T_{SB}^{2/1}$ be the binary subtree $T_{SB}$ rooted at the vertex $2/1$.  Since $u[]=1/2$ and $v[]=2/1$, applying induction to the correspondences above, it clearly follows that $$\left \{u[\bft],v[\bft]: \bft \mathrm{\hspace{1.5mm} has \hspace{1.5mm}length \hspace{1.5mm} of\hspace{1.5mm} }c \in \Z^+ \right \} = \left \{ f_{x_1} \cdots f_{x_c}(1/2), f_{x_1} \cdots f_{x_c}(2/1): x_i =L,R \right \}.$$  Thus, the generations of $T_{SB}$ can be taken to be equivalent to the generations of $X$ with respect to the second generation of $T_{SB}$ and the first generation of $X$; of note, however, is the fact that the parent of $m/n$ in $T_{SB}$ does not necessarily correspond to the parent of $[\min\{m,n\}, \max\{ m,n\}, m+n]$ in $X$.

\section{Application to the Three Hat Problem}
\setcounter{equation}{0}

%Change X from a player variable to Y to avoid confusion with $X$, $X$, and so on.

The final application of the Bernoulli process variation of the Fibonacci sequence that we explore is the ``Three Hat Problem'' puzzle.  In this puzzle, three players each have a positive integer on their respective hats and are told that two of the numbers are equal to the third.  Proceeding in a turn-wise, modular order, each player can either pass or announce his number if he has determined it explicity\footnote{It is generally assumed that each player is altruistic.}.  For further information on this puzzle, see \cite{BW,TR03, TR04, TR06}.

In \cite{TR03, TR04, TR06}, the Three Hat Problem puzzle is stated in the form \begin{equation}$$\text{``After} $n$ \text{rounds, player X concludes his number is} $m$\text{.  What are the other two numbers?''}$$  \end{equation} where the variables $n$, $m$, and X are varied over their respective parameters to increase or decrease the difficulty of the puzzle.  As an example of an application of the Bernoulli process variation of the Fibonacci sequence, we establish some  criteria for simplifying (5.1).

Note that, with respect to the three hat problem, we can express all non-base configurations as 

\noindent $\rho [a,b,c]$ for $a,b,c \in \Z^+$, $a+b=c$, and $\rho:\Z^3 \ra \Z^3$ a dimension permutation; in other words, we can each non-base configuration as a permuted state of $X$.  To give a more rigorous account of this fact, with respect to \cite{BW}, consider a three hat configuration $[a,b,c]$ with the restriction that $a \leq b \leq c$; we will refer to such a configuration as an {\em ordered configuration}.  Further, suppose that we reorder each configuration in the chain of the configuration $[a,b,c]$ so that each one is also an ordered configuration; we will refer to such a configuration chain as an {\em ordered configuration chain} or an {\em ordered chain} for brevity.  Since each non-base configuration chain contains the configuration $[1,2,3]$, we denote a configuration chain which omits the base configuration $[1,2,3]$ as {\em abbreviated}. 

Consider an ordered configuration chain as a linear graph where each configuration in the chain is represented by a labeled vertex; further, the configuration on the labels of the vertices are in the order of their location in the configuration chain and directed away from the base configuration $[1,1,2]$. Let $C_1, \ldots , C_n$ be ordered, abbreviated configuration chains; then let union of these chains, $\cup \{C_i : 1 \leq i \leq n \}$, be the graph on $\cup \{V(C_i): 1 \leq i \leq n \}$ where each $C_i$ is a subset of $\cup \{C_i : 1 \leq i \leq n \}$ and the indegree of each vertex of $\cup \{C_i : 1 \leq i \leq n \}$ is one.

\begin{prop}  The union of all primitive, abbreviated, ordered configuration chains is equal to $X$.\end{prop}%Redefine chain and these notions graph theoretically

\Proof Clearly, each vertex label ${\bf x}$ of $X$ is an element of $\Z^3$ in which the sum of the first two entries of ${\bf x}$ equal the third entry; thus, ${\bf x}$ is a three hat configuration.  In addition, by induction, since all pairwise combinations of entries of the first state vertex of $X$ have greatest common divisor equal to $1$, it follows that all pairwise combinations of entries of each state vertex along the path from the first state vertex up to ${\bf x}$ have greatest common divisor equal to $1$\footnote{Since, if $\gcd (a,b)=1$ and $a+b=c$, then $\gcd(a,c)=\gcd(b,c)=1$.}.  Therefore, ${\bf x}$ is a primitive, ordered configuration implying that $X$ is a subcollection of the union of all primitive, abbreviated, ordered configuration chains.

Now, each primitive, ordered three hat configuration ${\bf h}$ has an associated ordered chain which contains $[1,2,3]$.  Since the reduction operator generating this chain is the inverse of both transformations $\tau_0, \tau_1$ as described in proposition 2.2, it follows that the chain corresponding to ${\bf h}$ must be contained in $X$.\eproof

In \cite{BW}, a reduction scheme was given to reduce the all configurations to a base configuration.  This scheme involved a function\footnote{This function was denoted as $\sigma$ in \cite{BW}.} which, applied to a configuration $\bfs$, generates a new configuration where the largest entry of $\bfs$ is replaced by the absolute difference of the smaller two integer; this is, unless $\bfs$ is a base configuration in which case $\sigma(\bfs)=\bfs$.  The previous proposition allows us to view the tree of all possible primitive configurations of the ``Three Hat Problem'' puzzle in much the same way that we viewed $X$.  In other words, the reversal of the reduction scheme to solve problems of the form (5.1) is equivalent to considering all Bernoulli variations of the Fibonacci sequence herein.  Thus, the construction of $X$ in the previous section allows us to take advantage of the structure and results of the previous section to simply computations of solutions to (5.1). 

We say that a case must be {\em verified} when we must write out the dialogue for the players of a given configuration in order to tell if it is equivalent to (5.1).  As far as we know, a simple, deterministic method by which to solve for all configurations of (5.1) without also sifting through and checking a number of false configurations is currently unknown.  However, using the relation between the three hat problem and our Bernoulli process variation of the Fibonacci Sequence established in proposition 5.1, we develop criteria herein which allows us to exclude many configurations as possible solutions to (5.1) without checking them.  Each of the criterion introduced herein allows configurations to be excluded as solutions to (5.1) without being verified.  The reader should note that the criteria herein is preliminary and by no means do we make the claim that it is universally optimal for solving (5.1).

\begin{prop} The following configuration exclusion criteria are valid for (5.1).
\begin{enumerate}
	\item {\rm(Chain Length Bounds)} Let $L(\bfs)$ denote the length of the chain associated to a general configuration $\bfs$.  With respect to the number of rounds $n$, the following is valid based on X.
		\begin{itemize} 
		\item If player A is the solver, then all $\bfs$ such that $L(\bfs) < \left \lfloor \frac{3n-2}{2} \right \rfloor$ 			and all $\bfs$ such that $L(\bfs) > 3n-2$ can be excluded.  
		\item If player B is the solver, then all $\bfs$ such that $L(\bfs)<\left \lfloor \frac{3n-1}{2} \right \rfloor$ 				and all $\bfs$ such that $L(\bfs) > 3n-1$ can be excluded.
		\item If player C is the solver, then all $\bfs$ such that $L(\bfs)<\left \lfloor \frac{3n}{2} \right \rfloor$ and 			all $\bfs$ such that $L(\bfs) > 3n$ can be excluded.
		\end{itemize} 
	\item {\rm($m$ Lower Bound)} If a chain length $c$ satisfies criterion 1, but $c+2>m$, then all $\bfs$ such that $L(\bfs)=c$ can be excluded; similarly, if $\min \{L(\bfs): \bfs$ satisfies criterion $1\}+2>m$ satisfies criterion $1\}+3]<m$, then (5.1) has no solutions.
	\item {\rm ($m$-Prime Upper Bound)} Let $F[n]$ for integer $n$, be the $n$-th integer in the Fibonacci sequence.  If $m$ is prime and a chain length $d$ satisfies criterion 1, but $F[d+3]<m$, then all $\bfs$ such that $L(\bfs)=d$ can be excluded; similarly, if $F[\max\{L(\bfs): \bfs$ satisfies criterion $1\}+3]<m$, then (5.1) has no solutions for $m$ prime.
	\item {\rm ($m$-Equivalence Classes)} Group configurations into equivalence classes defined by $[\bfs]=\{\bfs':\max(\bfs)=\max(\bfs')\}$.  In other words, all configurations of one equivalence class have the exact same maximum entry.  Choose any and only one representative configuration from each class. 
\begin{itemize}
	\item Exclude all configurations in a class, if for that class, the representative configuration $\bfs$ is such that $\max(\bfs) > m$.
	\item Exclude all configurations in a class, if for that class, the representative configuration $\bfs$ is such that $\max(\bfs)$ does not divide $m$.
\end{itemize}
\item Check the two indices where player X is the largest integer in the remaining configurations not exahusted by criteria 1 through 4 to find all solutions of (5.1).
\end{enumerate}
\end{prop}

\begin{lem}
Consider statement (5.1).  For all general configurations that solve (5.1), if player A is the solver, then $$\left \lfloor \frac{3n-2}{2} \right \rfloor \leq L(\bfs) \leq 3n-2.$$  If player B is the solver, then $$\left \lfloor \frac{3n-1}{2} \right \rfloor \leq L(\bfs) \leq 3n-1.$$  If player C is the solver, then $$\left \lfloor \frac{3n}{2} \right \rfloor \leq L(\bfs) \leq 3n.$$
\end{lem}

We reference the reader to \cite{BW} for any terminology associated with the three hat problem which might not be clear with respect to the proof of the lemma.

\Proof We begin with player C.  To prove the lower bound, we consider the dialogue with the minimum number of cues within the given $n$ rounds.  Since each new cue is only dependent on the previous cue, we construct a dialogue which gives each cue to the last possible player while still following the optimal strategy outlined in \cite{BW}.  In round 1, the last possible person to have the cue is player C.  Then, since player C had the cue in round 1, he cannot have it again in round 2 implying that player A or player B must now have the cue.  Since player A has his turn before player B, in order to ensure the minimum number of cues, we assign player B the cue in round 2.  Following this format, we find that player A and player C will both have cues in round 3, player B will have a cue in round 4, and so on.
$$\begin{array}{clc}
&\text{Player A:} & \text{Pass}\\
\text{Round 1}&\text{Player B:} & \text{Pass}\\
&\text{Player C*:} & \text{Pass}\\
&\text{Player A:} & \text{Pass}\\
\text{Round 2}&\text{Player B*:} & \text{Pass}\\
&\text{Player C:} & \text{Pass}\\
&\text{Player A*:} & \text{Pass}\\
\text{Round 3}&\text{Player B:} & \text{Pass}\\
&\text{Player C*:} & \text{Pass}\\
$\vdots$&\hspace{8mm} $\vdots$ & $\vdots$\\
\end{array}$$From induction, we find that after the first round, player B has a cue in the even rounds and players B and C have a cue in the odd rounds.  Thus, overall, the minimum number of cues for odd $n$ rounds is $1+1+2+1+2+1+2+\cdots+1+2+1+1$\footnote{Since player C is the solver, we do not actually consider him to have a cue.} and the minimum number of cues for even $n$ rounds is $1+1+2+1+2+1+2+\cdots+1$.  Note that since player C is the solver, the game lasts the full $n$ rounds.

Now, for $n$ odd, noting that there are $(n-3)/2$ pairs of $1+2$ cues plus an additional $3$ cues, gives us $$3\left (\frac{n-3}{2}\right )+3=\frac{3(n-1)}{2}$$ possible cues. For $n$ even, there are $(n-2)/2$ pairs of $1+2$ cues plus an additional $2$ cues gives us $$3\left (\frac{n-2}{2}\right)+2=\frac{3n-2}{2}$$ cues.  Since the number of cues required to solve a configuration $\bfs$ is one less than $L(\bfs)$, we have that for $n$ odd, $$L(\bfs) \geq \frac{3n-1}{2}$$ and for $n$ even, $$L(\bfs) \geq \frac{3n}{2}.$$  Thus, we have proven the lower bounds.

To establish the upper bounds, we must maximize the number of cues for a given dialogue.  Since the only limitation on the structure of cues is that any one player cannot have two consecutive cues, clearly it is possible that each player has a cue in every possible round.  Thus, the dialogue becomes $$\begin{array}{clc}
&\text{Player A*:} & \text{Pass}\\
\text{Round 1}&\text{Player B*:} & \text{Pass}\\
&\text{Player C*:} & \text{Pass}\\
&\text{Player A*:} & \text{Pass}\\
\text{Round 2}&\text{Player B*:} & \text{Pass}\\
&\text{Player C*:} & \text{Pass}\\
&\text{Player A*:} & \text{Pass}\\
\text{Round 3}&\text{Player B*:} & \text{Pass}\\
&\text{Player C*:} & \text{Pass}\\
$\vdots$&\hspace{8mm} $\vdots$ & $\vdots$\\
\end{array}$$  Since player C is the solver, there are $3n-1$ cues implying that $L(\bfs) \leq 3n$ for both $n$ even and $n$ odd.

Note in statements where player A or player B is the solver in the $n$th round, the bounds change slightly due to the change in minimum and maximum cue structure.  So now, from the cases where player C was the solver, we can derive the bounds for the cases where player A or player B is the solver.

First, suppose that player A is the solver, then for $n$ odd, $$L(\bfs) \geq \frac{3n-1}{2}-1=\frac{3(n-1)}{2}$$ since player A normally has a cue in an odd round; for $n$ even, $$L(\bfs) \geq \frac{3n}{2}-1=\frac{3n-2}{2}$$  since player B does not get an opportunity to give a cue in the final round.  For the upper bound, since the game is shortened by two possible cues from the case where player C was the solver, $$L(\bfs) \leq 3n-2.$$  Using similar reasoning for player B as the solver, we have for $n$ odd, $$L(\bfs) \geq \frac{3n-1}{2}$$ and for $n$ even, $$L(\bfs) \geq \frac{3n-2}{2}.$$ Further, for an upper bound for player B as the solver, we obtain the bound $$L(\bfs) \leq 3n-1.$$\eproof

{\bf Proof of Proposition.} Criterion 1 is proven by the previous lemma.  Within this proof, we assume that $\bfs=S[\bft]$ for general edge traversal labeling $\bft$.

While we have denoted the chain length of a configuration $\bfs$ as $L(\bfs)$, we will denote the traversal length of a configuration as $\ell(\bfs)$.  Note that, in general, $L(\bfs)=\ell(\bfs)+1$.  Considering the problem by traversals, we can consider the bounds on $F[\alpha_1 \cdots \alpha_k]$ for a general integer $k$ by noting that, by definition, $F[\alpha_1 \cdots \alpha_{k-1} 0] < F[\alpha_1 \cdots \alpha_{k-1} 1]$.  Thus, the upper bound for an arbitrary traversal length $k$ is $$F[\underbrace{11\cdots 1}_{k}]$$ while the lower bound is $$F[\underbrace{00\cdots 0}_{k}].$$ Further, it is easy to show via induction that $$F[\underbrace{11\cdots 1}_{k}]=F_{k+4}$$ and $$F[\underbrace{00\cdots 0}_{k}]=k+3$$ where $F_{k+4}$ is the $(k+4)$-th integer in the Fibonacci sequence.  Now we proceed in proving criteria 2 and 3.

For criteria 2, $\min\{F[\bft]:\ell(\bft)=c\}=\ell(\bfs)+3=L(\bfs)+2$.  Thus, if $m<L(\bfs)+2$, then $m<F[\bft']$\footnote{So $F[\bft']$ cannot divide $m$.} for all $\bft'$ with $\ell(
S[\bft'])=c-1$.  The case where $\min \{L(\bfs): \bfs$ satisfies criterion $1\}+2>m$ tells us that this relation is true for all chain length $c$ such that configurations with chain length $c$ are not excluded by criterion 1, thus, by the first part of criterion 2, (5.1) would not have any solutions. For criterion 3, $\max\{F[\bft]:\ell(\bft)=d\}=F_{(d-1)+4}=F_{d+3}$; so if $m$ is prime and $m>F[\bft']$ for all $\bft'$ with $\ell(F[\bft'])=d-1$, then $F[\bft'] \neq m$.  Since $m$ is prime, $F[\bft']$ does not divide $m$ unless it is $m$, so all such $\bft'$ can be excluded.  The case where $F[\max\{L(\bfs): \bfs$ satisfies criterion $1\}+3]<m$ tells us that this relation is true for all chain lengths not excluded by criterion 1.  So, in this case, (5.1) has no solutions.

For $m$-Equivalence, it is easy to show that $[\bfs]$ is an equivalence class on $\HH$\footnote{The equality between maximum entries of primitive configurations is clearly reflexive, symmetric, and transitive.}.  Due to the primitivity assumption, it is a trivial observation that if  $\bfs$ solves (5.1), then $F[\bft]$ divides $m$\footnote{If $F[\bft]>m$, then $F[\bft]$ certainly cannot divide $m$}.  In the event that an indexing of $\bfs$ corresponds to (1.1) and $\max(\bfs)$ divides $m$, then it should be obvious that multiplying each entry of $\bfs$ by $m/\max(\bfs)$ gives an actual solution of (5.1).\eproof

While this proposition gives us theoretical criteria for excluding configurations as solutions to (5.1), we do not yet know enough about $X$ to take full advantage of all criteria in the proposition.  For instance, since the converse of the reflection theorem is not true, in general, finding a correspondence (if one exists) which links all equivalent values of a each code length would allow one to take full advantage of the $m$-equivalence classes; however, the reflection theorem can still be applied to the criteria.  As mentioned previously, we do not make the claim herein that our criteria is even the optimal criteria for reducing the complexity of solving (5.1).

\section*{Concluding Remarks}

Due to the introductory nature of this work, we do not take full advantage of any of the probablistic components of our construction, instead focusing on the structure of the collection of all variations of the Fibonacci sequence under this construction.  These probablistic components could be considered in future work.  Further, while we attempt to give some structural connections of our variation of the Fibonacci sequence, it goes without saying that there are more than likely connections or applications which we are unfamiliar with or have simply over-looked herein.

\noindent
{\bf Acknowledgement.} The author would like to thank Ian Fredenberg, Yang Wang, Ernie Croot, Prasad Tetali, Peter Winkler, and Donald Aucamp for very helpful discussions.  A large portion of the work herein was conducted while the author attended the REU program at the Georgia Institute of Technology in the summers of 2006 and 2007.

\end{document}